%regevHook.tex: Proof of a Conjecture of Amitai Regev About Hook multi-sets
%%a Plain TeX file by Doron Zeilberger (x pages)

%begin macros

\baselineskip=14pt
\parskip=10pt
\def\halmos{\hbox{\vrule height0.15cm width0.01cm\vbox{\hrule height
  0.01cm width0.2cm \vskip0.15cm \hrule height 0.01cm width0.2cm}\vrule
  height0.15cm width 0.01cm}}
\font\eightrm=cmr8 
\font\eighttt=cmtt8
\magnification=\magstephalf

\def\1{{\overline{1}}}
\def\2{{\overline{2}}}
\parindent=0pt
\overfullrule=0in

\def\frac#1#2{{#1 \over #2}}
%\headline={\rm  \ifodd\pageno  \RightHead  \else  \LeftHead  \fi}
%\def\RightHead{\centerline{
%Title
%}}
%\def\LeftHead{ \centerline{Doron Zeilberger}}
%end macros
\bf
\centerline
{
A Multi-set Identity for Partitions} \rm
\bigskip
\centerline{
{\it
Amitai REGEV\footnote{$^1$} {\eightrm \raggedright Department of
Theoretical Mathematics, Weizmann Institute of Science,
Rehovot 76100, ISRAEL. \hfill \break
amitai dot regev at weizmann dot ac dot il ,
http://www.wisdom.weizmann.ac.il/\~{}regev/ .}
and
Doron ZEILBERGER}\footnote{$^2$} {\eightrm \raggedright Department of
Mathematics, Rutgers University (New Brunswick), Hill Center-Busch
Campus, 110 Frelinghuysen Rd., Piscataway, NJ 08854-8019, USA.
{\eighttt zeilberg  at math dot rutgers dot edu} , \hfill \break
{\eighttt http://www.math.rutgers.edu/\~{}zeilberg/} .
Supported in part by the USA National Science Foundation. }
}

{\bf Added Sept. 22, 2009}

Guo-Niu Han kindly pointed out to us (something that  we should have noticed ourselves if
we would have been in the habit of reading {\bf carefully} all the papers
that we cite), that our main result is contained in [B.H.].

{\bf Introduction}

 Given an
integer-partition $\lambda \vdash n$ and a box (a cell) $v=[i,j]
\in \lambda$ it determines the arm length $a_v$ ($=\lambda_i-j$),
the leg length $l_v$ (=$\lambda'_j-i$), and the left length $f_v$
($=j-1$). Thus, for example, the hook length $h_v$ is given by
$h_v=a_v+l_v+1$. Denote $p_v=a_v+f_v+1$. C. Bessenrodt [B], and R.
Bacher and L. Manivel [B.M] (see also [B.H]) proved the following
identity:

$$(1)~~~~~~~~~~~~~~~~~~~~~~~~~~~~~~~~~~~~~~~~~~~~~~\sum_{\lambda\vdash
n}\sum_{v\in\lambda}x^{h_v}=\sum_{\lambda\vdash
n}\sum_{v\in\lambda}x^{p_v},~~~~~~~~~~~~~~~~~~~~~~~~~~~~~~~~~~~~~~$$
which is equivalent to the multi-set identity:
$$(2)~~~~~~~~~~~~~~~~~~~~~~~~~~~~~~~~~~~~~~\bigcup_{\lambda\vdash n}\{h_v\mid v\in\lambda\}=
\bigcup_{\lambda\vdash n}\{p_v\mid
v\in\lambda\}.~~~~~~~~~~~~~~~~~~~~~~~~~~~~~~~~~~$$ In this note we
prove the following refinement of (2).

Fill $v$ with a pair of numbers in two different ways:

{\bf First Filling}: Fill $v$ with $(a_v,l_v)$.

{\bf Second Filling}: Fill $v$ with $(a_v,f_v)$.

This yields the following two multi-sets of pairs:
$$
A_1(n)= \bigcup_{\lambda \vdash n} \{ ( a_v,l_v ) \vert v \in \lambda \} \quad ,
$$
$$
A_2(n)= \bigcup_{\lambda \vdash n} \{ ( a_v,f_v ) \vert v \in \lambda \} \quad .
$$

{\bf Theorem 1:} For all non-negative integers $n$ we have the
multi-set identity, $$A_1(n)=A_2(n).$$

\medskip
The proof here is by applying the technique of generating
functions. Theorem 1 indicates that for each $n$ there is a map
$\varphi$ on the cells of the partitions of $n$,
$\varphi:v\to\varphi(v)$, such that
$(a_v,f_v)=(a_{\varphi(v)},l_{\varphi(v)})$. The construction of
an explicit such $\varphi$ -- for all $n$ -- would yield a
bijective proof of Theorem 1.

\bigskip

{\bf The proof.}

\medskip

As usual, $(z)_a:=(1-z)(1-qz) \cdots (1-q^{a-1} z)$.

The proof would follow from the following two lemmas.

{\bf Lemma 1}: Let $M_1(c,d)(n)$ be the number of times the pair $(c,d)$ shows up in $A_1(n)$, then
$$
\sum_{n=0}^{\infty} M_1(c,d)(n)q^n = {{q^{c+d+1}} \over {1-q^{c+d+1}}} \cdot {{1} \over {(q)_\infty}} \quad.
\eqno(1)
$$

{\bf Lemma 2}: Let $M_2(c,d)(n)$ be the number of times the pair $(c,d)$ shows up in $A_2(n)$, then
$$
\sum_{n=0}^{\infty} M_2(c,d)(n)q^n = {{q^{c+d+1}} \over {1-q^{c+d+1}} } \cdot {{1} \over {(q)_\infty}} \quad.
\eqno(2)
$$

{\bf Proof of Lemma 2}: $M_2(c,d)(n)$ counts the number of Ferrers diagrams of $n$ where one of the cells
that has (right) arm $c$ and left-arm $d$ is {\bf marked}.
Obviously it belongs to a row of length $c+d+1$, and each such row
has exactly one such cell. Hence this is the same as counting the number of Ferrers diagrams of $n$ where one
of the rows of length $c+d+1$ is marked. We can construct such a Ferrers diagram (with any number of cells)
by first drawing that row of length $c+d+1$ (weight $q^{c+d+1}$) then putting {\bf below} it an arbitrary
Ferrers diagram with largest part $\leq c+d+1$, whose generating function is
$1/((1-q)(1-q^2) \cdots (1-q^{c+d+1}))$, and then placing {\bf above} the above-mentioned fixed row
any Ferrers diagram whose {\it smallest} part is $\geq c+d+1$,
whose generating function is
$1/((1-q^{c+d+1})(1-q^{c+d+2}) \cdots )$. Combining, we get that the
generating function
of such marked creatures, which is the left side of $(2)$, is the right side of $(2)$ \halmos.

Before proving Lemma 1 we have to recall certain basic facts from $q$-land.

{\bf Fact 1} (The $q$-Binomial Theorem [essentially Theorem 2.1  of
[A]\footnote{$^3$}
{\eightrm \raggedright But the ``conditions''  $|q|<1, |t|<1$, stated by Andrews, are, in our
world-view, a {\it category mistake}.} ]).
$$
{{1} \over {(z)_{a+1}}}=
\sum_{j=0}^{\infty} {{ (q)_{a+j}} \over {(q)_a (q)_j}} z^j \quad .
$$
(This is easily proved by induction on $a$).

When $a=\infty$ this simplifies to

{\bf Fact 2}
$$
{{1} \over {(z)_{\infty}}}=
\sum_{j=0}^{\infty} { {z^j} \over {(q)_j}} \quad .
$$

{\bf Fact 3}: The generating function for Ferrers diagrams bounded in an
$m$ by $n$ rectangle is ${ {(q)_{m+n}} \over {(q)_m (q)_n}}$.

This is Proposition 1.3.19 in [St] and Theorem 3.1 of [A]. Here is a proof by induction of
this elementary fact. Let the generating function be $F(m,n;q)$.
Consider the last cell of the top row. If it is occupied, the
generating function of these diagrams is $q^nF(m-1,n)$ (remove the
fully-occupied top row), if it is not, it is $F(m,n-1)$ (delete
the empty rightmost column), getting the recurrence
$F(m,n;q)=q^nF(m-1,n;q)+F(m,n-1;q)$. Then verify that the same
recurrence is satisfied by ${ {(q)_{m+n}} \over {(q)_m
(q)_n}}$, and check the trivial initial conditions $m=0$ and
$n=0$.

By sending $n$ to infinity we obtain

{\bf Fact 4}: The generating function for Ferrers diagrams with
parts bounded by $m$  is ${ 1 \over {(q)_m }}$.
By conjugation, this is also the
generating function for Ferrers diagrams with at most $m$ parts.

{\bf Proof of Lemma 1}: The left-side of $(1)$ is the generating function for Ferrers diagrams where one {\bf hook}
with arm-length $c$ and leg-length $d$ is {\bf marked}. Let's figure out the generating function (weight-enumerator)
for all such $(c,d)$-hook-marked Ferrers  diagrams.

  Suppose  the corner of that hook is at cell $(i+1,j+1)$ (i.e. the $(i+1)$-row and the $(j+1)$-column).
Here $0 \leq i < \infty$ and  $0 \leq j < \infty$.
Let's look at its {\it anatomy}. It consists of {\bf seven parts}.
(See diagram in   \hfill\break
http://www.math.rutgers.edu/\~{}zeilberg/mamarim/mamarimhtml/TemunaFerrers.html ).

{\bf 1.} Strictly {\bf left of} {\it and} {\bf above} cell
$(i+1,j+1)$. This is a fully occupied  $i$ by $j$ rectangle with
weight $q^{ij}$.

{\bf 2.} Above the arm (of length $c+1$). This is a fully occupied
$i$ by $c+1$ rectangle with weight $q^{(c+1)i}$.

{\bf 3.} To the left of the leg (of length $d+1$).
This is a fully occupied  $d+1$ by $j$ rectangle with weight $q^{(d+1)j}$ .

{\bf 4.} The Ferrers diagram with $\leq i$ rows lying {\bf above
and to the right} of the arm. By Fact 4, the generating function
of this is $1/(q)_i$.

{\bf 5.} The Ferrers diagram with $\leq j$ columns lying {\bf
below and to the left} of the leg. By Fact 4, the generating
function of this is $1/(q)_j$.

{\bf 6.} The hook itself. This gives generating function $q^{c+d+1}$.

{\bf 7.} The Ferrers diagram formed {\bf inside} the hook, i.e. lying below the
arm and to the right of the leg. By Fact 3 its generating function is
${ {(q)_{c+d}} \over {(q)_c (q)_d}}$.

Combining, we see that the generating function for these $(c,d)$-hook-marked Ferrers diagrams is
$$
{ {(q)_{c+d}} \over {(q)_c (q)_d} } \cdot q^{c+d+1} \cdot q^{ij+i(c+1)+j(d+1)} \cdot {{1} \over {(q)_i}}
\cdot {{1} \over {(q)_j}} \quad .
$$

Summing over {\bf all} $0 \leq i,j < \infty$, we get that the generating function on the left of $(1)$ equals
$$
{ {(q)_{c+d}} \over {(q)_c (q)_d} }  q^{c+d+1}
\sum_{i=0}^{\infty} \sum_{j=0}^{\infty}
 q^{ij+i(c+1)+j(d+1)} {{1} \over {(q)_i}}
 {{1} \over {(q)_j}}
$$
$$
=
{ {(q)_{c+d}} \over {(q)_c (q)_d} } q^{c+d+1}
\sum_{i=0}^{\infty}
{{1} \over {(q)_i}}
q^{i(c+1)}
\sum_{j=0}^{\infty}
q^{j(i+d+1)}
{{1} \over {(q)_j}} \quad
$$
$$
=
{ {(q)_{c+d}} \over {(q)_c (q)_d} } q^{c+d+1}
\sum_{i=0}^{\infty}
{{1} \over {(q)_i}}
q^{i(c+1)} {{1} \over {(q^{d+i+1})_\infty}} \quad ,
$$
by Fact 2 with $z=q^{d+i+1}$. This, in turn, equals
$$
{ {(q)_{c+d}} \over {(q)_c } } q^{c+d+1} \sum_{i=0}^{\infty}
q^{i(c+1)} {{1} \over {(q)_{\infty}}} {{(q^{i+1})_{d}} \over
{(q)_d} }
$$
$$
=
{{q^{c+d+1}} \over {(q)_\infty}}
{ {(q)_{c+d}} \over {(q)_c } } \sum_{i=0}^{\infty}  q^{i(c+1)} {{(q)_{i+d}} \over {(q)_d (q)_i} }
$$
$$
={{q^{c+d+1}} \over {(q)_\infty}}
{ {(q)_{c+d}} \over {(q)_c } }    {{1} \over {(q^{c+1})_{d+1} } } \quad ,
$$
by Fact 1 with $z=q^{c+1}$. Finally, this equals
$$
={{q^{c+d+1}} \over {(q)_\infty}} \cdot { {(1-q)(1-q^2) \cdots (1-q^{c+d}) } \over {(1-q)(1-q^2) \cdots (1-q^{c+d+1})} }
={{q^{c+d+1}} \over {(q)_\infty}} \cdot {{1} \over {(1-q^{c+d+1})} } \quad \halmos .
$$

{\bf Reference}

[A] G. Andrews, {\it The Theory of Partitions,} Cambridge
University Press (1984).

[B] C. Bessenrodt, {\it On hooks of Young diagrams,} Ann. of
Comb., {\bf 2} (1998), pp. 103-110.

[B.M]  R. Bacher and L. Manivel, {\it Hooks and powers of parts in
partitions}, Sem. Lothar. Combin. Vol {\bf 47}, article B47d,
(2001) 11 pages.

[B.H] C. Bessenrodt and G.N Han, {\it Symmetry distribution
between hook length and part length for partitions}, Discrete
Mathematics (online 4 June 2009).

[St] R. Stanley, {\it Enumerative Combinatorics I}, Cambridge
Studies in Advanced Mathematics 49 (1986).

\end